\documentclass[11pt]{article}
\usepackage{amsfonts,amssymb,graphicx}
\usepackage{epsfig,amsmath,latexsym}
\usepackage{amscd, amssymb, amsthm, amsmath,eucal, verbatim}
\usepackage{epstopdf}
\newtheorem{thm}{Theorem}[section]

\begin{document}

\newcommand{\G}{{\mathbb G}}
\newcommand{\Z}{{\mathbb Z}}
\newcommand{\R}{{\mathbb R}}
\newcommand{\I}{{\mathrm{Id}}}
\newcommand{\lk}{{\mathcal{L}}}
\newcommand{\sL}{{\mathcal{L}}}
\newcommand{\sC}{\mathcal{C}}
\newcommand{\sD}{{\mathcal{D}}}
\newcommand{\sG}{{\mathcal{G}}}
\newcommand{\sE}{{\mathcal{E}}}
\newcommand{\sH}{{\mathcal{H}}}
\newcommand{\dg}{{\mathcal{D}}}

\newcommand{\ct}{\cite}
\newcommand{\pr}{\protect\ref}
\newcommand{\su}{\subseteq}
\newcommand{\pa}{{\partial}}

\newcounter{numb}
\date{\today}
 
\title{Unknot diagrams requiring a quadratic number of 
Reidemeister moves to untangle}
\author{Joel Hass\footnote {Supported in part by NSF grant DMS 3289292.}  ~and Tahl Nowik}

\date{November 14, 2007}
 
\maketitle

\begin{abstract}
We present a sequence of diagrams of the unknot for which the minimum number of Reidemeister moves required to pass to the trivial diagram
is quadratic with respect to the number of crossings. These bounds apply both in $S^2$ and in $\R^2$.
\end{abstract}

\maketitle

\section{Introduction}\label{A}

In this paper we give a family of unknot diagrams $D_n$ with  $D_n$ having $7n-1$ crossings and with at  least $ 2n^2 + 3n - 2$ Reidemeister moves required to transform $D_n$ to the trivial diagram. 
These are the first examples for which a non-linear lower bound has been established. 

A knot in $\R^3$ or $S^3$ is commonly represented by a knot diagram, a generic projection of
the knot to a plane or 2-sphere. A diagram is an
immersed oriented planar or spherical curve with finitely many double points, called crossings.
Each crossing is marked to indicate a strand, called the overcrossing, that lies
above the second strand, called the undercrossing. The original knot can be recovered, up to isotopy, by constructing a curve with the overcrossing arcs pushed 
slightly above the plane of the diagram and the remainder of the diagram lying in this plane.

\begin{figure}[htbp]
\centering
\scalebox{0.4}{\includegraphics{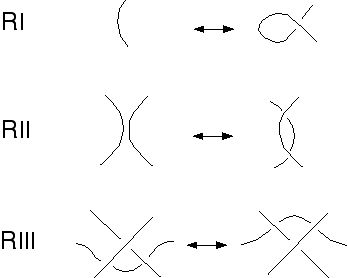}}
\caption{The three types of Reidemeister moves }
\label{Reidmoves}
\end{figure}

Alexander and Briggs~\cite{AleBri26} and independently Reidemeister~\cite{R1} showed that two diagrams 
of the same knot can be connected through a sequence of moves of three types, commonly called Reidemeister moves, shown in Figure~\ref{Reidmoves}.
The number of such moves required to connect two equivalent diagrams is difficult to  estimate. 
An exponential upper bound is obtained in~\cite{HasLagPip}, where it is shown that there is a positive constant $c$ 
such that given an unknot diagram $D$ with $n$ crossings, no more than $2^{c n}$ Reidemeister moves are required  to transform $D$ 
to the trivial knot diagram.

We can get some lower bounds by looking at classical invariants of diagrams such as crossing numbers, writhes and  winding numbers, since a single Reidemeister move changes these numbers by 0, 1 or 2. However non-trivial lower bounds are difficult to obtain.  This is a common situation in complexity theory. While upper bounds can be established by careful analysis of one procedure, lower bounds require
somehow bounding from below the running time of all possible procedures. It is quite difficult to get lower bounds even for a  particular pair of equivalent diagrams, as seen in  \cite{Carter}, \cite{Hagge}, \cite{Hayashi} and \cite{Ostlund}. We note that examples are constructed in \cite{HassSnoeyinkThurston} that show that it may require exponentially many faces to construct a PL spanning disk for an unknotted polygon. However these examples can be transformed to the trivial diagram using only a linear number of Reidemeister moves.
  
Given two knot diagrams $D,E$ of the same knot, we define the \emph{Reidemeister distance} 
$d(D, E)$ between $D$ and $E$ to be the minimal number of Reidemeister moves required to pass from $D$ to $E$.
One may consider this notion in either $S^2$ or $\R^2$, and our result will hold in both settings.
Our main tool is an invariant of knot diagrams developed in \cite{HassNowik}, 
and used there to obtain new linear lower bounds on the  Reidemeister distance.
 
\section{The diagrams}

\begin{figure}[htbp]
\centering
\scalebox{0.6}{\includegraphics{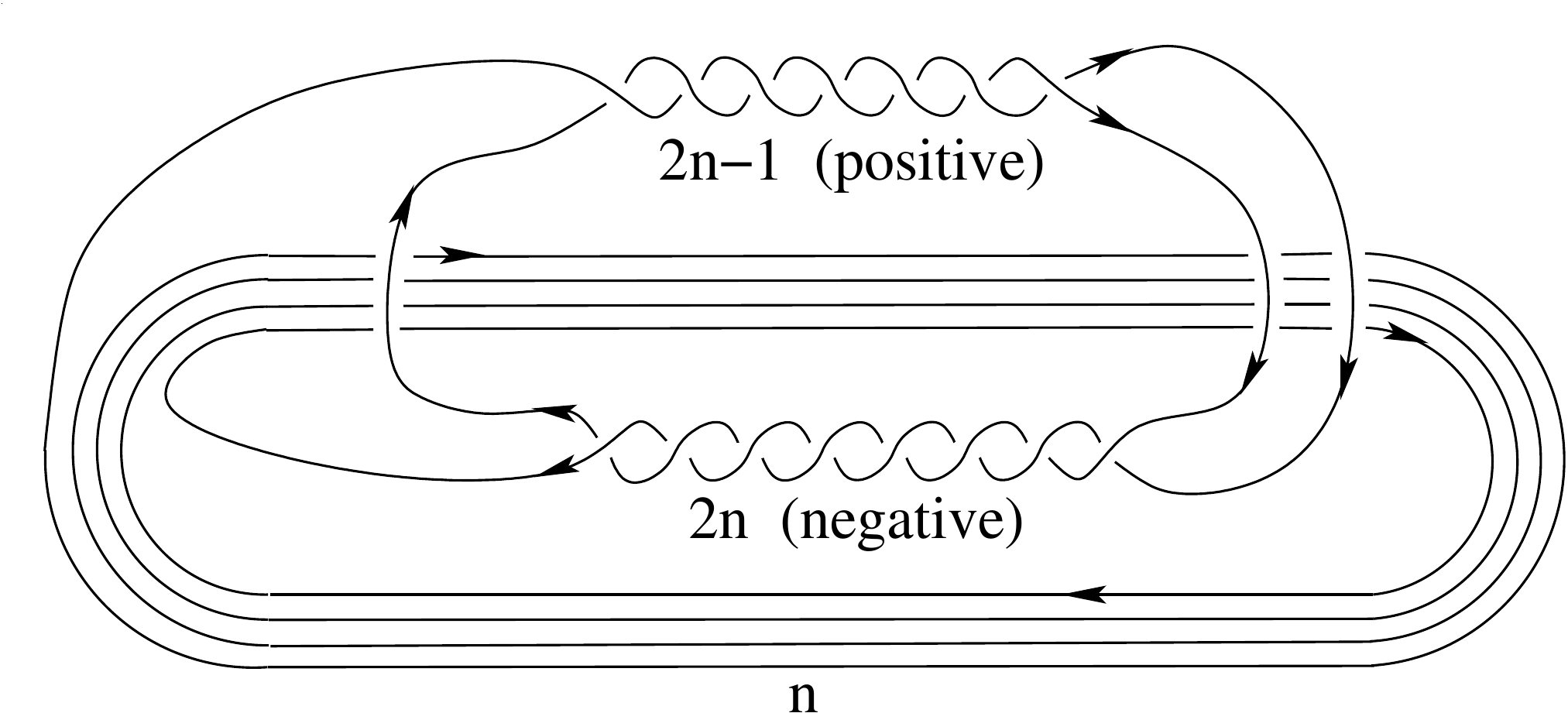}}
\caption{The diagram $D_n$ for $n=4$.}
\label{f1}
\end{figure}

Let $U$ denote the trivial knot diagram.
We will present a sequence $D_n$ of diagrams of the unknot, for which $d(D_n , U)$,
the number of Reidemeister moves required to pass from  $D_n$ to $U$
grows quadratically with respect to the number of crossings of $D_n$.
More precisely, we  prove: 

\begin{thm}\label{main}
In both $S^2$ and $\R^2$, the diagram $D_n$ of Figure~\ref{f1}, which has $7n-1$ crossings, satisfies:
$$ 2n^2 + 3n - 2 \ \ \  \leq \ \ \  d(D_n, U) \ \ \  \leq \ \ \  2n^2 + 3n.$$
\end{thm}

We will prove the lower bound in $S^2$, and the upper bound in $\R^2$, and it follows
that both bounds hold in both settings.
We recall the definition of the invariant of knot diagrams in $S^2$ introduced in \ct{HassNowik}.
We denote the set of all knot diagrams in $S^2$ by $\dg$ and
the set of all two component links in $\R^3$ by $\lk$. 
Given a knot diagram $D \in \dg$ and a crossing $a$ in $D$, define $D^a \in \lk$
to be the two component link obtained by smoothing the crossing $a$.
Given a knot diagram $D$, let   $D_+$ denote the set of all positive
crossings in $D$ and  $D_-$ the set of all negative crossings.
Given an invariant of two component links $\phi:\lk \to S$ where $S$ is any set,
let $\G_S$ be the free abelian group with basis $\{X_s, Y_s\}_{s \in S}$. We then define 
the invariant $I_\phi : \dg \to \G_S$ to be
$$I_\phi(D) = \sum_{a \in D_+} X_{\phi(D^a)}  + \sum_{a \in D_-} Y_{\phi(D^a)}.$$

In this work $\phi$ is taken to be the linking number, $lk:\lk\to\Z$,  giving the invariant
$I_{lk} : \dg \to \G_\Z$. 
In \ct{HassNowik} it is shown that the change in the value of $I_{lk}$ 
resulting from a Reidemeister move has one of the following 
forms:
\begin{enumerate}
\item For an RI move: $X_0$ or $Y_0$. 
\item For an RII move: $X_k + Y_k$ or $X_k + Y_{k+1}$.
\item For an RIII move: $X_k - X_{k+1}$ or $Y_k - Y_{k+1}$.  
\end{enumerate}
Let $R$ be the set of elements in $\G_\Z$ of the
form $X_0$, $Y_0$, $X_k+Y_k$, $X_k+Y_{k+1}$, $X_k - X_{k+1}$, $Y_k - Y_{k+1}$, and their negatives. 
That is, $R$ is the set of all elements of $\G_\Z$ that may appear as the change in the value of 
$I_{lk}(D)$ as the result of performing a Reidemeister move on $D$. The set $R$ generates $\G_\Z$ and
the length of an element of $\G_\Z$ with respect to this generating 
set is called its {\em $R$-length}. Given two diagrams $D, E$ of the same knot, the $R$-length
of $I_{lk}(D) - I_{lk}(E)$ is a lower bound for $d(D,E)$ in $S^2$, and therefore also in $\R^2$.
In particular, if $D$ is a diagram of the unknot, then since $I_{lk}(U)=0$,
the $R$-length of $I_{lk}(D)$ gives a lower bound for $d(D,U)$. 
We use this procedure to obtain our lower bound for $d(D_n,U)$.

\noindent
{\bf Proof of Theorem~\ref{main}.}
A computation shows that 
$$I_{lk}(D_n)=nX_n + nX_{-n} + (2n-1)X_{-1} + 3nY_0.$$
Indeed, each crossing in 
the top horizontal string of $2n-1$ crossings contributes $X_{-1}$, 
each crossing in the bottom horizontal string of $2n$ crossings contributes $Y_0$,
each crossing in the left vertical line of $n$ crossings contributes $Y_0$,  
each crossing in the middle vertical 
line of $n$ crossing contributes $X_n$, and
each crossing in the right vertical line of $n$ crossings contributes $X_{-n}$.
Together this gives $nX_n + nX_{-n} + (2n-1)X_{-1} + 3nY_0$, and we
denote this element of $\G_\Z$ by $v_n$.

We prove Theorem \pr{main} by first showing that the $R$-length of $v_n$ is
at least $2n^2 + 3n - 2$, and then demonstrating an explicit sequence of $2n^2 + 3n$
Reidemeister moves in $\R^2$ from $D_n$ to $U$.
 
Let $g:\G_\Z \to\Z$ be the homomorphism defined by
$g(X_k)=1+|k|$ and $g(Y_k)=-1-|k|$ for all $k$.
Then $|g(r)| \leq 1$ for all $r \in R$,
and $g(v_n)=2n^2 + 3n -2$. It follows that the $R$-length  of $v_n$ is 
at least $2n^2 + 3n - 2$.

It remains to demonstrate an explicit sequence of $2n^2 + 3n$ Reidemeister moves in $\R^2$, from
$D_n$ to the trivial diagram. 
Start by sliding the top
horizontal string of $2n-1$ crossings in Figure~\ref{f1} in the clockwise direction, 
across the $n$ horizontal strands. This requires $n(2n-1)$ RIII moves.
Then cancel these $2n-1$ positive crossings with $2n-1$ of the negative 
crossings, now lying to the left of them,
via $2n-1$ RII moves, arriving at the diagram in Figure \pr{f2}.
With $n$ additional RII moves we arrive at the diagram in Figure \pr{f3}. 
Finally,  perform $n+1$ RI moves to get to the trivial diagram, resulting in a total of $2n^2 + 3n$ moves.

\begin{figure}[htbp]
\centering
\scalebox{0.3}{\includegraphics{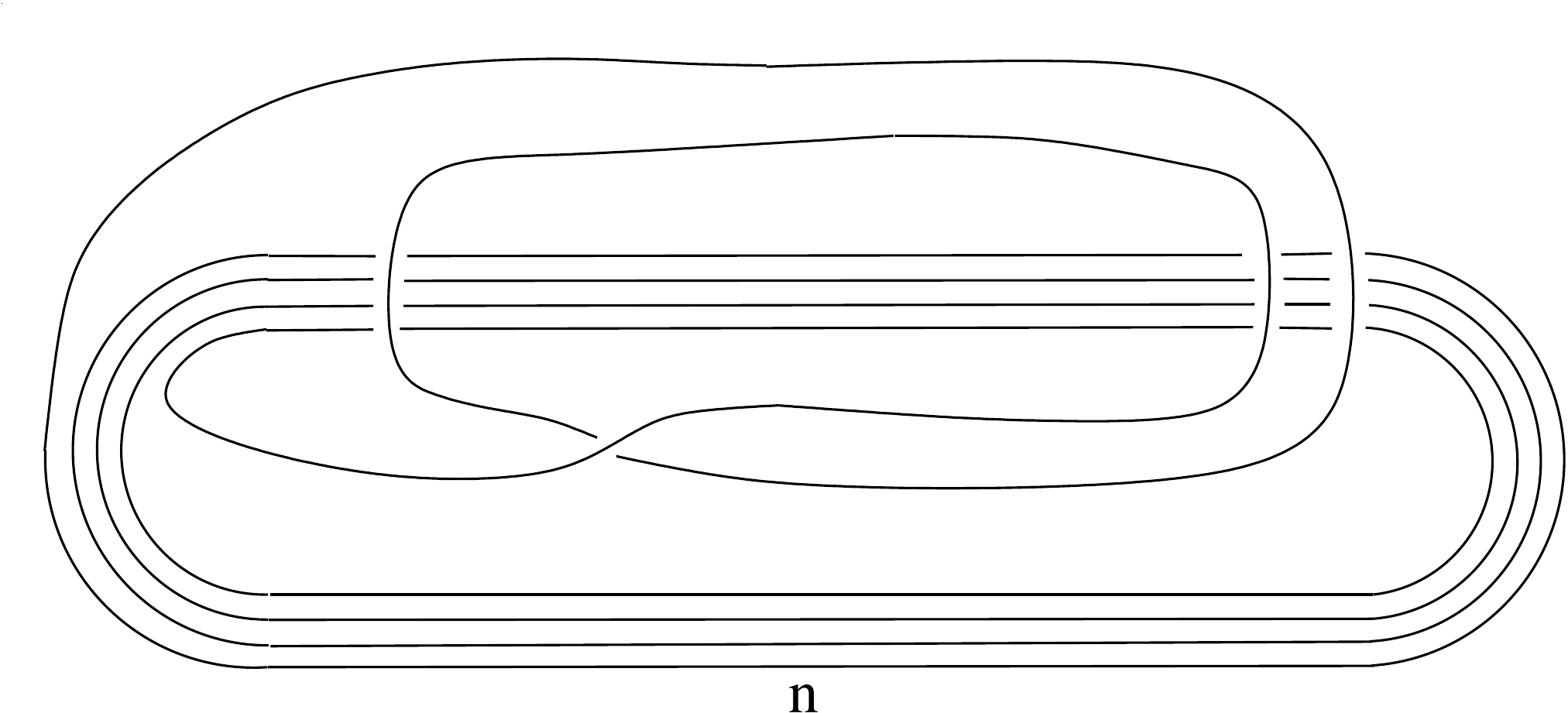}}
\caption{$D_n$ after $2n^2-n$ RIII moves and $2n-1$ RII moves.}
\label{f2}
\end{figure}

\begin{figure}[htbp]
\centering
\scalebox{0.3}{\includegraphics{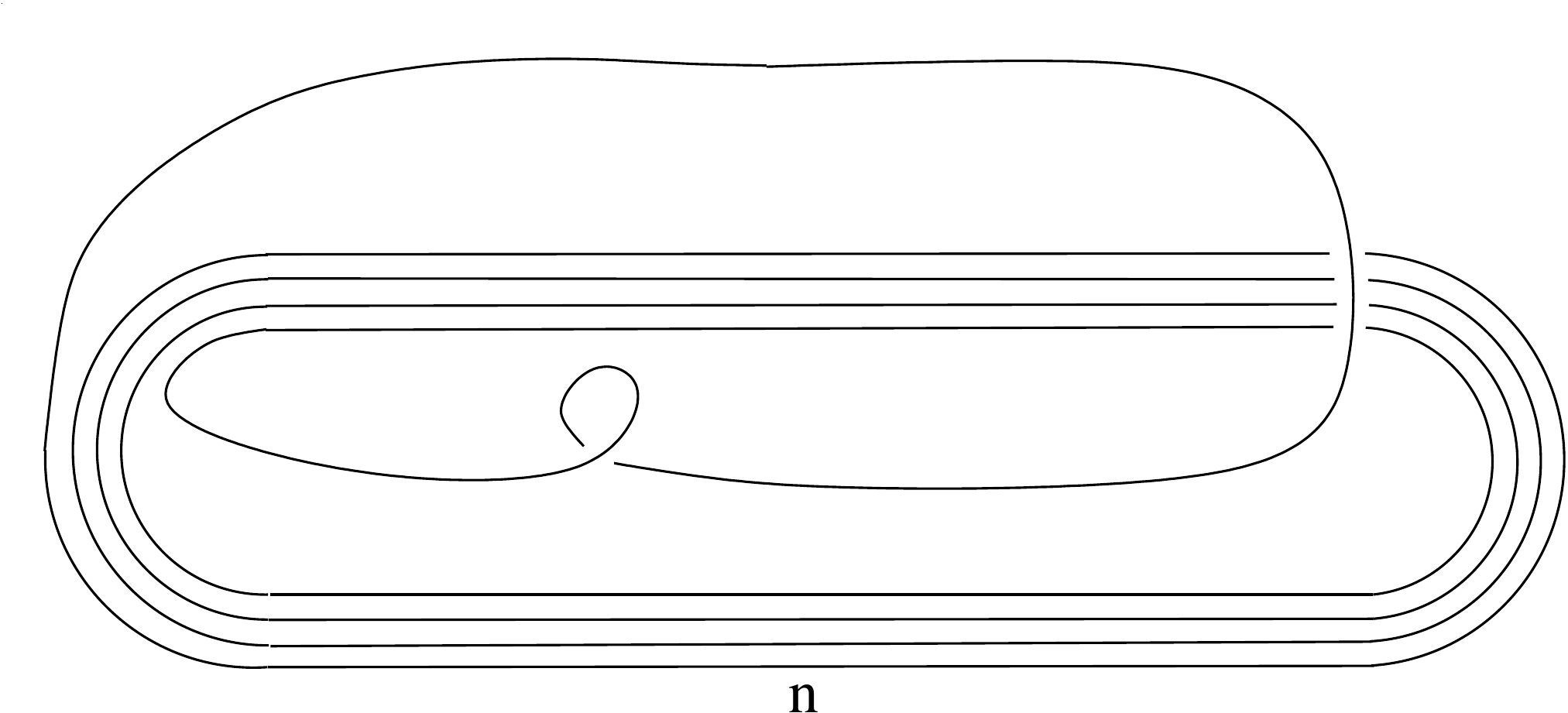}}
\caption{$D_n$ after additional $n$ RII moves.}
\label{f3}
\end{figure}

\qed

\vspace*{.5in}
\begin{tabular}{ll}
{\tt email}: & {\tt hass@math.ucdavis.edu} \\ & {\tt tahl@math.biu.ac.il } \\
~~~ \\
~~~ \\
{\tt address}: & J. Hass \\
& Department of Mathematics \\
& University of California, Davis \\
& Davis, CA 95616 \\
~~~ \\
& T. Nowik \\
& Department of Mathematics \\
& Bar-Ilan University \\
& Ramat-Gan 52900, Israel \\
 \\
\end{tabular}


\begin{thebibliography}{StoPC}

\bibitem{Adams}
C. Adams,
{\em The Knot Book. An elementary introduction to the mathematical theory
of knots},
W. H. Freeman, New York, 1994.

\bibitem{AleBri26}
J. W. Alexander and G. B. Briggs,
{\em On types of knotted curves},
Ann. Math., {\bf 28} (1926/27), 562--586.

\bibitem{Arnold94}
V. I. Arnold: 
{\em Plane Curves, Their Invariants, Perestroikas and Classifications,}
 Advances in Soviet Mathematics, {\bf  21} (1994), 33-91.

\bibitem{Carter}
J. S. Carter, M. Elhamdadi ,  M. Saito and S. Satoh,
{\em A lower bound for the number of Reidemeister moves of type III}
Topology and its Applications, {\bf 153} (2006) 2788--2794.

\bibitem{Hagge} 
T.J. Hagge, 
{\em Every Reidemeister move is needed for each knot type,}
arxives:math.GT/0404145.

\bibitem{HasLagPip}
J. Hass and J. C. Lagarias,
{\em The number of Reidemeister moves needed for unknotting,}
J. Amer. Math. Soc. {\bf 14} (2001), no. 2, 399--428.

\bibitem{HassNowik}
J. Hass, T. Nowik: ``Invariants of knot diagrams'' 
arXiv:0708.2509

\bibitem{HassSnoeyinkThurston}
J. Hass, J. Snoeyink and W.P. Thurston
{\em The Size of Spanning Disks for Polygonal Curves,} 
Discrete \& Computational Geometry {\bf 29} (2003) (1): 1-17. 

\bibitem{Hayashi}
C. Hayashi,  
{\em A lower bound for the number of Reidemeister moves for unknotting.} 
J. Knot Theory Ramifications {\bf 15} (2006), no. 3, 313--325.

\bibitem{Ostlund}
O. \"{O}stlund,
{\em Invariants of knot diagrams and relations among Reidemeister
moves,} 
J. Knot Theory Ramifications, {\bf 10} (2001), no. 8, 1215--1227. 

\bibitem{R1}
H. Reidemeister, 
{\em Knoten und Gruppen}, Abh. Math. Sem., Univ. Hamburg,
{\bf 5} (1926), 7--23.

\end{thebibliography}
\end{document}